\documentclass[11pt]{article}
\usepackage{latexsym,amssymb,amsmath,amscd,amsthm,amsxtra}
\usepackage{mathrsfs}
\usepackage{epsfig}
\usepackage{pgf,tikz}
\usepackage{graphicx}
\usepackage{graphics}
\usepackage{setspace}
\usepackage[all,cmtip]{xy}
\newcommand{\counte}{section}
\newtheorem{define}{\bf Definition}[\counte]

\newtheorem{prop}{\bf Proposition}[\counte]
\newtheorem{lemma}{\bf Lemma}[\counte]
\newtheorem{theorem}{\bf Theorem}[\counte]
\newtheorem{coro}{\bf Corollary}[\counte]
\newtheorem{example}{\bf Example}[\counte]

\def\r{{\mathrm{rank}}}
\def\Z{{\mathbb Z}}
\def\Q{{\mathbb Q}}
\def\R{{\mathbb R}}

\newcommand{\D}{\displaystyle}

\author{Zhao Xu-an, zhaoxa@bnu.edu.cn\\ Ruan Yangyang, 201831130006@mail.bnu.edu.cn\\ Wang Ran, 201621130041@mail.bnu.edu.cn\\
School of Mathematical Sciences, Beijing Normal University\\Key Laboratory
of Mathematics and Complex Systems\\ Ministry of Education,
China, Beijing 100875}

\title{The isomorphism types of parabolic subgroups in Kac-Moody groups\thanks{The authors are supported by National Science Foundation of China, 11571038.}}

\date{}

\begin{document}

\maketitle
\begin{abstract}
In this paper we discuss the isomorphism types of parabolic subgroups in Kac-Moody groups. The results have applications in the study of topology of Kac-Moody groups and their classifying spaces.
\end{abstract}



\section{Introduction}

Let $A=(a_{ij})$ be an $n\times n$ integer matrix which satisfies

(1) $a_{ii}=2$, for $1\leq i\leq n$;

(2) $a_{ij}\leq 0$, for $1\leq i\not = j\leq n$;;

(3) $a_{ij}=0$ if and only if $a_{ji}=0$,

\noindent then $A$ is called a Cartan matrix.

For a Cartan matrix $A$, let $h$ be the complex vector space with basis $\Pi^{\vee}=\{\alpha^{\vee}_1,\alpha^{\vee}_2,\cdots,\alpha^{\vee}_n\}$. We denote the dual basis of $\Pi^\vee$ in the dual vector space $h^*$ by $\{\omega_1,\omega_2,\cdots,\omega_n\}$, i.e. $\omega_i(\alpha_j^{\vee})=\delta_{ij}$ for $1\leq i,j\leq n$. Let $\Pi=\{\alpha_1,\cdots,\alpha_n\}\subset h^*$ be given by $\langle\alpha^{\vee}_i,\alpha_j\rangle=a_{ij}$ for all $i,j$, then $\alpha_i=\sum\limits_{j=1}^n a_{ji}\omega_j$. When the Cartan matrix $A$ is singular, then $\{\alpha_i,1\leq i\leq n\}$ is not a basis of $h^*$. $\alpha_i, \alpha_i^\vee, \omega_i,1\leq i\leq n$ are called respectively the simple roots, simple coroots and fundamental dominant weights.

By the work of Kac\cite{Kac_68} and Moody\cite{Moody_68}, it is well known that for each Cartan matrix $A$, there exists a complex Lie algebra $g(A)$ which is called Kac-Moody Lie algebra.


The Kac-Moody Lie algebra $g(A)$ is generated by $\alpha^\vee_i,e_i,f_i,1\leq i\leq n$ over the complex field $\mathbb{C}$,
with the defining relations:

(1) $[\alpha_i^{\vee},\alpha_j^{\vee}]=0$ for all $i,j$;

(2) $[e_i,f_j]=\delta_{ij}\alpha_i^{\vee}$ for all $i,j$;

(3) $[\alpha_i^{\vee},e_j]=a_{ij}e_j,[\alpha_i^{\vee},f_j]=-a_{ij}f_j$ for all $i,j$;

(4) $\mathrm{ad}(e_i)^{-a_{ij}+1}(e_j)=0$, for all $1\leq i\not=j\leq n$;

(5) $\mathrm{ad}(f_i)^{-a_{ij}+1}(f_j)=0$, for all $1\leq i\not=j\leq n$.

Kac and Peterson\cite{Kac_Peterson_83}\cite{Kac_Peterson_84}\cite{Kac_85} constructed the simply connected Kac-Moody group $G(A)$ with Lie algebra $g(A)$.

Cartan matrices and their associated Kac-Moody Lie algebras, Kac-Moody groups are divided into three types.

(1) Finite type, if $A$ is positive definite. In this case, $G(A)$ is just the simply connected complex semisimple Lie group with Cartan matrix $A$.

(2) Affine type, if $A$ is positive semi-definite and has rank $n-1$.

(3) Indefinite type otherwise.



Let $S$ be the set $\{1,2,\cdots,n\}$. For each proper subset $I\subsetneq S$, the matrix $A_I=(a_{ij})_{i,j\in I}$ is also a Cartan matrix. And there is a Kac-Moody subalgebra $g(A_I)$ of $g(A)$ which is generated by $\alpha_i^{\vee},e_i,f_i,i\in I$, and a parabolic subalgebra $g_I(A)$ generated by $e_i,i\in I$ and $\alpha_i^{\vee},f_i,1\leq i\leq n$. Corresponding to $g(A_I)$ and $g_I(A)$, there is a connected subgroup $G(A_I)$ and a parabolic subgroup $G_I(A)$ in $G(A)$.


There is a canonical anti-linear involution on the Lie algebra $g(A)$ which is defined by $\sigma(e_i)=f_i,\sigma(f_i)=e_i$ and $\sigma(\alpha_i^\vee)=-\alpha_i^\vee$. The fixed real subalgebra $k(A)$ is spanned by $\sqrt{-1}\alpha_i^\vee, 1\leq i\leq n; e_{\alpha}+f_{\alpha}, \sqrt{-1}(e_{\alpha}-f_{\alpha}),\alpha\in \Delta_+$(the positive root system of $G(A)$). $k(A)$ has a real Cartan subalgebra $t$ which is generated by $\sqrt{-1}\alpha_i^\vee,1\leq i\leq n$. And $\sigma$ induces an automorphism $w$ on $G(A)$. The unitary form $K(A)$ of $G(A)$ is defined as the fixed subgroup $G(A)^w$ and $k(A)$ is its Lie algebra. The unitary form of the Kac-Moody subgroup $G(A_I)$ is $K(A_I)=K(A)\cap G(A_I)$ with Lie algebras $k(A_I)=k(A)\cap g(A_I)$. And the unitary form of the parabolic subgroup $G_I(A)$ is $K_I(A)=K(A)\cap G_I(A)$ with Lie algebras $k_I(A)=k(A)\cap g_I(A)$. For $I=\emptyset$, $K_\emptyset (A)$ is the maximal torus subgroup $T$ of $K(A)$ and $t$ is its Lie algebra. The groups $K(A_I)$ and $K_I(A)$ have the same simple root system $\Pi_I=\{\alpha_i|i\in I\}$ and Dynkin diagram. The readers may refer to \cite{Kac_Peterson_83}\cite{Kac_Peterson_84}\cite{Kac_85}\cite{Kumar_02} for details.



Let $t_I=k(A_I)\cap t$ be the Cartan subalgebra of $K(A_I)$ and $t'_I$ be the central subalgebra of $k_I(A)$. Let $T_I$ and $T'_I$ be the torus subgroups of $K(A)$ corresponding to $t_I$ and $t'_I$. Since $K(A_I)$ commutes with $T'_I$, we have a homomorphism $\pi_I: K(A_I)\times T'_I \to K_I(A),(a,b)\mapsto ab$. And $\pi_I$ induces a Lie algebra homomorphism $d\pi_I: k(A_I)\times t'_I \to k_I(A),(X,Y)\mapsto X+Y$.

Generally, $d\pi_I$ is neither injective nor surjective. However, we will prove that $d\pi_I$ is an isomorphism if and only if the Cartan matrix $A_I$ is non-singular. In this case $d\pi_I$ is surjective and the kernel of $\pi_I$ is a discrete subgroup $\Gamma_I$ in the center of $K(A_I)\times T'_I$, thus $K_I(A)\cong K(A_I)\times T'_I/\Gamma_I$. Then a natural problem is

\noindent{\bf Problem: }For the pair $(A,I)$ with $\det A_I\not=0$, how to determine the discrete subgroup $\Gamma_I$.

If $\det A_I=0$, $K_I(A)$ can't be written as a quotient of $K(A_I)\times T'_I$, but as a semi-direct product of $K(A_I)$ and $T_{S-I}$.

For Kac-Moody groups of finite type, Duan and Liu\cite{Duan_Liu_12} determine the isomorphism types of centralizer of an element in a simply connected compact Lie group. Their results contain the computation of the isomorphism type of parabolic subgroups of Kac-Moody group of finite type. In this paper, we will generalize their methods to general cases.

The methods to determine the isomorphism type of $K_I(A)$ is as follows. First we observe that the kernel of $\pi_I$ equals to the kernel of its restricted homomorphism on maximal torus. Then by identifying the Cartan subalgebra $t_I\oplus t'_I$ of $k(A_I)\times t'_I$ with the Cartan subalgebra $t$ of $k_I(A)$, and comparing the unit lattice $L(A)\subset t_I\oplus t'_I$ with the unit lattice $L_I(A)\subset t$, we determine the group $\Gamma_I$. We have

\noindent {\bf Theorem 1: }
Let $L(A)$ and $L_I(A)$ be the unit lattices of $T$ and $T_I\times T'_I$ respectively. If $\det A\not=0$, $\Gamma_I$ is isomorphic to $L(A)/L_I(A)$.

\noindent {\bf Theorem 2: }
Let $E_I(A)$ be the subgroup of $\mathbb{Z}^{|I|}$ generated by the row vectors $e_1,e_2,\cdots,e_n$ of matrix $\widetilde A_{I}=(a_{ij})_{i\in S,j\in I}$ and $E(A_I)$ be the subgroup generated by the row vectors $e_{i_1},e_{i_2},\cdots,e_{i_{|I|}}$ of matrix $A_{I}=(a_{ij})_{i,j\in I}$, then we have $\Gamma_I\cong E_I(A)/E(A_I)$.

The contents of this paper are as follows. In section 2, we discuss the computation of the center of a Kac-Moody group. In section 3, we determine the group $\Gamma_I$. 





\section{The center of the Kac-Moody group $K(A)$}

In this section we compute the center $Z(K(A))$ of the Kac-Moody group $K(A)$.

By definition $$Z(K(A))=\{a\in K(A)|aba^{-1}=b,\forall b\in K(A)\}.$$ By the general theory of Kac-Moody groups, the center of $K(A)$ is contained in the maximal torus subgroup $T$. Hence an element $a\in Z(K(A))$ can be written as $a=\exp(\Lambda)$ with $\Lambda\in t$. And $$Z(K(A))=\{a|Ad(a)Y=Y,\forall Y\in k(A)\}.$$ Therefore $a\in Z(K(A))$ if and only if $Ad(a)e_i=e_i$ and $Ad(a)f_i=f_i$ for all $i$. And we have $Ad(a)e_i=Ad(\exp (\Lambda))e_i=\exp(ad(\Lambda))e_i=e^{\alpha_i(\Lambda)}e_i$. Similarly $Ad(a)f_i=e^{-\alpha_i(\Lambda)}f_i$. Hence $a\in Z(K(A))$ if and only if $\alpha_i(\Lambda)\in 2\pi \sqrt{-1} \Z$.

From now on we may omit the constant $2\pi \sqrt{-1}$ temporarily for convenience and we can recover it at any time. An element
$\Lambda=\sum\limits_{i=1}^n \lambda_i \alpha_i^{\vee}\in t$ can be regarded as a row vector $(\lambda_1,\lambda_2,\cdots,\lambda_n)$ under the basis $\Pi^{\vee}$. Similarly the root $\alpha_i=\sum\limits_{j=1}^n a_{ji} w_j $ can be regarded
as a column vector $(a_{1i},a_{2i},\cdots,a_{ni})^T$, which is exact the $i$-th column vector of $A$.
\begin{define}
For a Cartan matrix $A$, we define $\mathcal{Z}(A)=\{\Lambda\in t|\alpha_i(\Lambda)\in \Z, 1\leq i\leq n\}$ and $L(A)=\{\Lambda\in t|\exp(\Lambda)=e\}$. $L(A)$ is the kernel of the exponential map $\exp: t\to T$ and is called the unit lattice of $T$.
\end{define}

It is easy to see that $L(A)\subset \mathcal{Z}(A)$. And for the simply connected Kac-Moody group $K(A)$, $L(A)$ is spanned by $\alpha_1^{\vee},\alpha_2^{\vee},\cdots,\alpha_n^{\vee}$.

\begin{prop}
The center of Kac-Moody group $K(A)$ is isomorphic to $\mathcal{Z}(A)/L(A)$.
\end{prop}

\noindent{\bf Proof: } We consider the surjective homomorphism $\phi: \mathcal{Z}(A)\to Z(G(A)),\Lambda\mapsto \exp(\Lambda)$. It is easy to see that the kernel of $\phi$ is $L(A)$. This proves the proposition.

The condition $\alpha_i(\Lambda)\in \Z, 1\leq i\leq n$ is equivalent to equation $(\lambda_1,\lambda_2,\cdots,\lambda_n)A=(k_1,k_2,\cdots,k_n)$ with $k_1,k_2,\cdots,k_n\in \Z$. Let $K=(k_1,k_2,\cdots,k_n)$, the equation can be written as
\begin{equation}  \label{equation1}
\Lambda A=K
\end{equation}
If $\det A\not=0$, then the solution of Equation \ref{equation1} is $\Lambda=KA^{-1}$. It can be written as $\Lambda=\sum\limits_{i=1}^n {k_i} v_i$ with $v_i$ being the $i$-th row vector of matrix $A^{-1}$. Hence $v_i,1\leq i\leq n$ generates $\mathcal{Z}(A)$. Let $\bar v_i$ be the image of $v_i$ under the homomorphism $\phi$, then $\bar v_i,1\leq i\leq n$ generates the group $Z(K(A))$. \qed

Let $t_{\Q}$ be the $\Q$-vector space $L(A)\otimes \Q$. We define a linear map $R: t_{\Q}\to \Q^n,  \Lambda\mapsto \Lambda A$.
Let $e_1,e_2,\cdots,e_n$ be the row vectors of $A$. They span a lattice subgroup $E$ in $\Z^n$. we have

\begin{coro}
If $\det A\not=0$, then $Z(K(A))\cong \Z^n/E$ and it is a finite group with order $|\det A|$.
\end{coro}

\noindent {\bf Proof:} If $\det A\not=0$, then $R$ is an isomorphism. Hence $\mathcal{Z}(A)/L(A)\cong R(\mathcal{Z}(A))/R(L(A))=\Z^n/E$. \qed

If $\det A=0$, Let $F$ be the intersection of $\Z^n\subset \Q^n$ with the $\Q$-vector space spanned by $e_1,e_2,\cdots,e_n$, we have

\begin{coro}
If $\det A=0$, then $Z(K(A))$ is isomorphic to the product of a finite group $Z'(K(A))$ and a torus subgroup $T^{n-k}\subset T$, where $Z'(K(A))\cong F/E$ and $k=rank(A_I)$.\label{coro}
\end{coro}

\noindent{\bf Proof: }We have two exact sequences.
$$0\to\ker R \to \mathcal{Z}(A)\to \mathrm{im} R\cap \Z^n\to 0,\ \  0\to\ker R\cap L(A)\to L(A)\to \mathrm{im} L(A)\to 0,\ \  $$ The quotient of these two sequences gives an exact sequence $$0\to T^{n-k}\to Z(K(A))\to F/E\to 0$$ Since $T^{n-k}$ is divisible, the exact sequence is split. This proves the corollary.\qed


\begin{example}
Let $A=(a_{ij})_{2\times 2}$ be a Cartan matrix of rank $2$, we compute the center of $K(A)$.

$\Lambda=(\lambda_1,\lambda_2)\in \mathcal{Z}(K(A))$ if and only if
$$\D{\left\{
  \begin{array}{ll}
   2\lambda_1+a_{21}\lambda_2=k_1 \\
   a_{12}\lambda_1+2\lambda_2=k_2
  \end{array}
\right.}$$
If $\Delta=4-a_{12}a_{21}\not=0$, then
$$\D{\left\{
  \begin{array}{ll}
   \lambda_1=\frac{2}{\Delta}k_1+\frac{-a_{21}}{\Delta}k_2 \\
   \lambda_2=\frac{-a_{12}}{\Delta}k_1+\frac{2}{\Delta}k_2
  \end{array}
\right.}$$
Hence we have $\D{(\lambda_1,\lambda_2)=k_1(\frac{2}{\Delta},\frac{-a_{12}}{\Delta})+k_2(\frac{-a_{21}}{\Delta},\frac{2}{\Delta})}$.

We set $\D{v_1=(\frac{2}{\Delta},\frac{-a_{12}}{\Delta})}$ and $\D{v_2=(\frac{-a_{21}}{\Delta},\frac{2}{\Delta})}$, then $\bar v_1$ and $\bar v_2$ are the generators of group $Z(K(A))$. Explicit computation shows

If $a_{12}$ is odd, then $\bar v_1$ is the generator of $Z(K(A))$ with order $| \Delta |$. And $Z(K(A))$ is isomorphic to $Z_{| \Delta |}$.

If $a_{21}$ is odd, then $\bar v_2$ is the generator of $Z(K(A))$ with order $| \Delta |$. And $Z(K(A))$ is isomorphic to $Z_{| \Delta |}$.

If both $a_{12}$ and $a_{21}$ are even, then $\bar v_1$ and $\bar v_2+\frac{a_{21}}{2}\bar v_1$ are the generators of $Z(K(A))$ with order $\frac{| \Delta |}{2}$ and order $2$. And $Z(K(A))$ is isomorphic to $Z_{\frac{| \Delta |}{2}}\times Z_{2}$.

If $\Delta=0$, then $a_{12}a_{21}=4$ and we must have $A=\left(
                                                                     \begin{array}{cc}
                                                                       2 & -2 \\
                                                                       -2 & 2 \\
                                                                     \end{array}
                                                                   \right)$ or $\left(
                                                                     \begin{array}{cc}
                                                                       2 & -1 \\
                                                                       -4 & 2 \\
                                                                     \end{array}
                                                                   \right)$




For $A=\left(
         \begin{array}{cc}
           2 & -2 \\
           -2 & 2 \\
         \end{array}
       \right)$, $e_1=(2,-2),e_{2}=(-2,2)$ and $E=\Z(2,-2),F=\Z(1,-1)$, hence $Z'(K(A))=\Z_2$.

For $A=\left(
         \begin{array}{cc}
           2 & -1 \\
           -4 & 2 \\
         \end{array}
       \right)$, $e_1=(2,-1),e_{2}=(-4,2)$ and $E=\Z(2,-1),F=\Z(2,-1)$, hence $Z'(K(A))=\{1\}$.

       In both cases $k=1$.
\end{example}

An $n\times n$ integer matrix $A$ has a decomposition $A=P\widetilde A Q$, where $P,\widetilde A,Q$ are integer matrices with $\det P,\det Q\in \{1,-1\}$ and $\widetilde A=\mathrm{diag}[q_1,q_2,\cdots,q_n]$ with $q_1\mid q_2\mid\cdots \mid q_n$. If $\det A=0$ and $\r (A)=k$, then $q_{k+1}=q_{k+2}=\cdots=q_n=0$. If $A$ is a Cartan matrix, then Equation \ref{equation1} becomes $\Lambda P\widetilde A Q=K$. Denoting $\Lambda P$ and $KQ^{-1}$ by $\widetilde \Lambda$ and $\widetilde K$ respectively, then we have $\widetilde \Lambda \widetilde A=\widetilde K$. An easy check shows that

\begin{prop}
$Z(K(A))\cong Z_{q_1}\times Z_{q_2}\times \cdots\times Z_{q_k}\times T^{n-k}$. 
\end{prop}

This provides an alternative proof for Corollary \ref{coro}

\section{The kernel $\Gamma_I$ of the homomorphism $\pi_I$}
First we compute the center $Z_I(A)$ of the parabolic subgroup $K_I(A)$. For $I\subsetneq S$, let $\mathcal{Z}_I(A)=\{\Lambda\in t|\alpha_i(\Lambda)\in \Z, \forall i\in I\}$.

\begin{lemma}
The center $Z_I(A)$ of $K_I(A)$ is $\exp( \mathcal{Z}_I(A))$.
\end{lemma}

\noindent{\bf Proof: }Since $Z_I(A)\subset T$, an element $a\in Z_I(A)$ can be written as $a=\exp(\Lambda)$. A similar argument as previous section shows that $a\in Z_I(A)$ if and only if $\Lambda\in \mathcal{Z}_I(A)$. This gives a surjective homomorphism $\phi:\mathcal{Z}_I(A)\to Z_I(A), \Lambda\mapsto \exp(\Lambda)$ and proves the lemma.\qed

Let $A$ be a Cartan matrix and $I\subsetneq S$ be a proper subset with $\det A_I\not=0$.
We have defined the group homomorphism $\pi_I: K(A_I)\times T'_I\to K_I(A), (a,b)\mapsto ab$. It induces a Lie algebra homomorphism $d\pi_I: k(A_I)\oplus t'_I\to k_I(A),(X,Y)\mapsto X+Y$. By checking the definition we have


\begin{lemma}
1. The Lie algebra $t'_I=\{\Lambda\in t|\alpha_i(\Lambda)=0, \forall i\in I\}$, or equivalently $\Lambda=(\lambda_1,\lambda_2,\cdots,\lambda_n)$ $\in t'_I$ if and only if $\Lambda \widetilde A_I=0$, where $\widetilde A_I=(a_{ij})_{i\in S,j\in I}$.

2. The kernel of $d\pi_I$ is $\{(X,-X)|X\in t_I\cap t'_I\}$ and it is isomorphic to $t_I\cap t'_I$.

3. The kernel of $\pi_I$ is $\{(a,a^{-1})|a\in T_I\cap T'_I\}$ and it is isomorphic to $T_I\cap T'_I=K(A_I)\cap T'_I$.

\end{lemma}

\begin{lemma}
$d\pi_I$ is an isomorphism if and only if $\det A_I\not=0$
\end{lemma}

\noindent{\bf Proof: }Since $k_I(A)$ and $k(A_I)$ have the same Dynkin diagram, by the previous lemma, it is sufficient to prove that $t_I\cap t'_I=\{0\}$ and $t_I+t'_I=t$.







We set $l=\r(A_I),m=\r(\widetilde{A}_I)$, then $l\leq m\leq |I|$. By 1 and 2 of Lemma 3.2, the dimensions of $t'_I$ and $t_I\cap t'_I$ are $n-m$ and $|I|-l$. Therefore

$\dim (t_I+ t'_I)=\dim t_I+\dim t'_I-\dim t_I\cap t'_I=|I|+n-m-(|I|-l)=n-(m-l)$.

Hence $d\pi_I$ is an isomorphism if and only if $|I|-l=0$ and $n-(m-l)=n$. So we have $rank(A_I)=l=m=|I|$. This means that $\det A_I\not=0$. \qed


In the following we need only consider the case when $\det A_I\not=0$.

Let $I=\{i_1,i_2,\cdots,i_{|I|}\},i_1<i_2<\cdots<i_{|I|}$ and $J=S-I=\{j_1,j_2,\cdots,j_{n-|I|}\},j_1<j_2<\cdots<j_{n-|I|}$, then the column vectors $\alpha_{i_1},\alpha_{i_2},\cdots \alpha_{i_{|I|}}$ form the $n\times |I|$ matrix $\widetilde A_I$. The matrix $\widetilde A_I$ can be decomposed into two sub-matrices. The $i_1,i_2,\cdots,i_{|I|}$-th rows of $\widetilde A_I$ form the matrix $A_I$, and the $j_1,j_2,\cdots,j_{n-|I|}$-th rows form another matrix $\bar A_I$. By the decomposition $t=t_I\oplus t_J$, each $\Lambda\in t$ can be written as $(\Lambda_1,\Lambda_2)$ with $\Lambda_1\in t_I$ and $\Lambda_2\in t_J$.
By definition if $\Lambda\in \mathcal{Z}_I(A)$, then $\Lambda \widetilde A_I=(k_{i_1},k_{i_2},\cdots,k_{i_{|I|}}), k_{i_1},k_{i_2},\cdots,k_{i_{|I|}}\in \Z$. Let $K_I$ be the row vector $(k_{i_1},k_{i_2},\cdots,k_{i_{|I|}})$, then we can write the equation $\Lambda \widetilde {A_I}=K_I$ as
\begin{equation}
\Lambda_1 A_I+\Lambda_2 \bar A_I=K_I.
\end{equation}\label{equation2}

The vector $\Lambda$ can also be decomposed as $\Lambda=\Lambda_1^I+\Lambda_2^I$ with $\Lambda_1^I\in t_I$ and $\Lambda_2^I\in t'_I$. The relation for the two decompositions are described in the following lemma.

\begin{lemma}
For $\Lambda=(\Lambda_1,\Lambda_2)$, Let $K_I=\Lambda \widetilde{A}_I=\Lambda_1 A_I+\Lambda_2 \bar A_I$. If $\det A_I\not=0$, then $\Lambda_1^{I}=(K_I A_I^{-1},0)$ and $\Lambda_2^{I}=(-\Lambda_2\bar A_I A_I^{-1},\Lambda_2)$.
\end{lemma}

\noindent{\bf Proof: }We need to show that $\Lambda_1^{I}\in t_I$, $\Lambda_2^{I}\in t'_I$ and $\Lambda=\Lambda_1^I+\Lambda_2^I$.
$\Lambda_1^{I}=(K_I A_I^{-1},0)\in t_I$ is obvious. $\Lambda_2^{I}\in t'_I$ is given by $\Lambda^I_2 \widetilde A_I=-\Lambda_2\bar A_I A_I^{-1} A_I+\Lambda_2 \bar A_I=0$. By $\Lambda_1 A_I+\Lambda_2 \bar A_I=K_I$ we have $\Lambda_1+ \Lambda_2 \bar A_I A_I^{-1}=K_I A_I^{-1}$. Hence $\Lambda=(\Lambda_1,\Lambda_2)=(K_I A_I^{-1}-\Lambda_2\bar A_I A_I^{-1},\Lambda_2)=\Lambda_1^I+ \Lambda_2^I$. This proves the lemma. \qed

We have the following commutative diagram

\begin{center}
$\xymatrix{
    t_I\oplus t'_I  \ar[d]_{id} \ar[r]^{\exp_I} & T_I\times T'_I  \ar[d]_{\bar\pi_I} \ar[r]^{\subset} & K(A_I)\times T'_I \ar[d]^{\pi_I} \\
    t \ar[r]^{\exp} & T \ar[r]^{\subset} & K_I(A)}$
\end{center}

\noindent where $\exp_I=(\exp,\exp):t_I\oplus t'_I\to T_I\times T'_I$.
The homomorphism $\pi_I$ is restricted to a homomorphism $\bar \pi_I: T_I\times T'_I\to T$ between maximal tori.

It is obvious that $\ker \pi_I=\ker \bar\pi_I$. The element $a=\exp_I(\Lambda)$ lies in $\ker \bar \pi_I$ if and only if $\exp(\Lambda)=e$, i.e. $\Lambda$ lies in $L(A)$.

\begin{define}
For $I\subset S$, we define $L_I(A)=\{\Lambda=\Lambda_1^I+\Lambda_2^I \in t_I\times t'_I|\exp(\Lambda_1^I)=e,\exp(\Lambda_2^I)=e\}$, which is called the unit lattice of $T_I\times T'_I$.
\end{define}

The lattice $L_I(A)$ is isomorphic to the lattice $(L(A)\cap t_I)\oplus (L(A)\cap t'_I)$. By Lemma 3.4 and the definition of $L_I(A)$, we have a simple but useful result
\begin{coro}
For an element $\Lambda=(\lambda_1,\lambda_2,\cdots,\lambda_n)\in L(A)$, $K_I=\Lambda \widetilde{A}_I$, then $\Lambda \in L_I(A)$ if and only if $K_I A_I^{-1}\in \Z^{|I|}$ or $\Lambda_2 \bar{A}_I A_I^{-1}\in \Z^{n-|I|}$.
\end{coro}

Now we have

\noindent {\bf Theorem 1: }
Let $L(A)$ and $L_I(A)$ be the unit lattices of $T$ and $T_I\times T'_I$. If $\det A_I \not=0$, then $\Gamma_I=\exp_I(L(A))\cong L(A)/L_I(A)$.



\noindent{\bf Proof: }Since the exponential map $\exp_I$ is surjective, each element $a\in \Gamma_I$ can be written as $a=\exp_I(\Lambda)$. By considering the left square of the previous commutative diagram we have $\exp(\Lambda)=\bar \pi_I(\exp_I(\Lambda))=\bar \pi_I(a)=\{e\}$. So $\Lambda\in L(A)$. As a
result, we can define a surjective homomorphism $f: L(A)\to \Gamma_I, \Lambda\mapsto \exp_I(\Lambda)$. It is obvious that the kernel of $f$ is $L_I(A)$. This proves the theorem. \qed

\noindent {\bf Theorem 2: }
Let $E_I(A)$ be the subgroup of $\mathbb{Z}^{|I|}$ generated by the row vectors $e_1,e_2,\cdots,e_n$ of the matrix $\widetilde A_{I}=(a_{ij})_{i\in S,j\in I}$ and $E(A_I)$ be the subgroup generated by the row vectors $e_{i_1},e_{i_2},\cdots,e_{i_{|I|}}$ of the matrix $A_{I}=(a_{ij})_{i,j\in I}$, then we have $\Gamma_I\cong E_I(A)/E(A_I)$.

\noindent{\bf Proof: }Let $\phi: \Z^n\to \Z^{|I|}$ be the $\Z$-linear map given by $\Lambda\mapsto \Lambda\widetilde A_I$, then $\phi(L(A))$ is a subgroup of $\Z^{|I|}$ generated by $(1,0,\cdots,0)\widetilde A_I,(0,1,\cdots,0)\widetilde A_I,\cdots,(0,0,\cdots,1)\widetilde A_I$, which are just the $n$ vectors $e_1,e_2,\cdots,e_n$. For an element $\Lambda\in \Z^n$, let $K_I=\Lambda \widetilde A_I$.
Then by Corollary 3.1, $\Lambda\in L_I(A)$ if and only if $K_I A_I^{-1}\in Z^{|I|}$. This shows that $K_I$ is in the image of map $\Z^n\to Z^{|I|},\Lambda\mapsto \Lambda \left(
                                                                           \begin{array}{c}
                                                                             A_I \\
                                                                             O \\
                                                                           \end{array}
                                                                         \right)$. So $\phi(L_I(A))$ is $E_I(A)$.

For $\Lambda\in \ker \phi$, we have $\Lambda \widetilde A_I=0$, that is $\Lambda_1A_I+\Lambda_2 \bar A_I=0$. This shows that $(\Lambda_2 \bar A_I A_I^{-1},0)=(-\Lambda_1,0)\in L(A)$. Therefore $\Lambda \in L_I(A)$. So $\ker \phi\subset L_I(A)$. Hence we have $L(A)/L_I(A)\cong L(A)/\ker \phi\Big{/}L_I(A)/\ker \phi \cong E_I(A)/E(A_I)$. \qed



Let $V$ be an $n$-dimensional real vector space with a lattice subgroup $L$. We assume $V$ has a direct sum decomposition $V=V_1\oplus V_2$, then we can define an abelian group $\Gamma(V,L;V_1,V_2)$ which is given by $L/(V_1 \cap L\oplus V_2\cap L)$.
The group $\Gamma(V,L;V_1,V_2)$ is trivial if and only if $L=V_1 \cap L\oplus V_2\cap L$, i.e. the direct sum decomposition of $V$ induces a direct sum decomposition of $L$. By this definition we have $\Gamma_I=\Gamma(t,L(A),t_I,t'_I)$.



\begin{example}
Let $V=\R^n$, and $L$ be a lattice generated by standard basis $e_1,e_2,\cdots, e_n$ of $V$. If $V_1$ is spanned by $e_1-e_2,e_2-e_3,\cdots,e_{n-1}-e_{n}$ and $V_2$ is spanned by $e_1+e_2+\cdots+e_n$, then $\Gamma(V,L,V_1,V_2)\cong \Z_n$. And the natural homomorphism $\pi:V_1/(L\cap V_1)\times V_2/(L\cap V_2)\to V/L$ is an n-fold cyclic covering.
\end{example}

This example corresponds to the fact that $U(n)\cong SU(n)\times S^1/\Z_n$.

In the following we compute $\Gamma_I$ for $3\times 3$ Cartan matrix.

\begin{example}
Let $A=(a_{ij})$ be a $3\times 3$ Cartan matrix, we compute the isomorphism type of the parabolic subgroup $K_I(A)$. By the symmetry of indices $1,2,3$, we only need to compute the cases $I=\{1\}$ and $\{1,2\}$.

For $I=\{1\}$, the groups $\mathcal{Z}_I(A)=\{\Lambda=(\lambda_1,\lambda_2,\lambda_3)|a_{11}\lambda_1+a_{21}\lambda_2+a_{31}\lambda_3=k_1\in\mathbb{Z}\}$.
$t_I$ is generated by $\alpha_1^\vee$(or $(1,0,0)$) and $t'_I=\{\Lambda=(\lambda_1,\lambda_2,\lambda_3)|a_{11}\lambda_1+a_{21}\lambda_2+a_{31}\lambda_3=0\}. $

By Corollary 3.1 the element $\Lambda=(\lambda_1,\lambda_2,\lambda_3)\in L_I(A)$ if and only if $\lambda_1,\lambda_2,\lambda_3\in \mathbb{Z}$ and $a_{21}\lambda_2+a_{31}\lambda_3\in 2\mathbb{Z}$, if and only if

$$\left\{
  \begin{array}{ll}
    \lambda_{1}, \lambda_{2},\lambda_3\in \Z & \hbox{if both $a_{21}$ and $a_{31}$ are even;} \\
    \lambda_{1}, \lambda_{2}\in\Z ,\lambda_3\in 2\Z, & \hbox{if $a_{21}$ is even and $a_{31}$ is odd;} \\
    \lambda_{1}, \lambda_3\in \Z, \lambda_{2}\in 2\Z,  & \hbox{if $a_{21}$ is odd and $a_{31}$ is even;} \\
    \lambda_{1}, \lambda_{2},\lambda_3\in \Z, \lambda_{2}+\lambda_3\in 2\Z, & \hbox{if both $a_{21}$ and $a_{31}$ are odd.}
  \end{array}
\right.$$

By $a_{11}\lambda_1+a_{21}\lambda_2+a_{31}\lambda_3=k_1$, we get $\lambda_1=\frac{k_1}{2}+\frac{-a_{21}}{2}\lambda_2+\frac{-a_{31}}{2}\lambda_3$.

Hence $(\lambda_1,\lambda_2,\lambda_3)=k_1(\frac{1}{2},0,0)+\lambda_2(-\frac{a_{21}}{2},1,0)+\lambda_3(-\frac{a_{31}}{2},0,1)=k_1(\frac{1}{2},0,0)+(-\frac{\lambda_2 a_{21}+\lambda_3 a_{31}}{2},1,0)$. This shows that $t'_I$ is spanned by $(-\frac{a_{21}}{2},1,0)$ and $(-\frac{a_{31}}{2},0,1)$.

If both $a_{21}$ and $a_{31}$ are even, then $K_I(A)$ is isomorphic to $SU(2)\times T^2$. Otherwise $K_I(A)$ is isomorphic to $SU(2)\times T^2/\mathbb{Z}_2$.


For $I=\{1,2\}$, $(\lambda_1,\lambda_2,\lambda_3)\in \mathcal{Z}_I(A)$ satisfies
$$\left\{
  \begin{array}{ll}
   a_{11}\lambda_1+a_{21}\lambda_2+a_{31}\lambda_3=k_1\in \Z \\
   a_{12}\lambda_1+a_{22}\lambda_2+a_{32}\lambda_3=k_2\in \Z
  \end{array}
\right.$$

We set $\Delta=\left|
             \begin{array}{cc}
               a_{11} & a_{12} \\
               a_{21} & a_{22} \\
             \end{array}
           \right|\not=0$, $\Delta_1=\left|
             \begin{array}{cc}
               a_{11} & a_{12} \\
               a_{31} & a_{32} \\
             \end{array}
           \right|$ ,$\Delta_2=\left|
             \begin{array}{cc}
               a_{21} & a_{22} \\
               a_{31} & a_{32} \\
             \end{array}
           \right|$ ,then
$$\left\{
  \begin{array}{ll}
   \lambda_1=\frac{a_{22}}{\Delta}(k_1-a_{31}\lambda_3)+\frac{-a_{21}}{\Delta}(k_2-a_{32}\lambda_3) \\
   \lambda_2=\frac{-a_{12}}{\Delta}(k_1-a_{31}\lambda_3)+\frac{a_{11}}{\Delta}(k_2-a_{32}\lambda_3)
  \end{array}
\right.$$

Hence we have

$(\lambda_1,\lambda_2,\lambda_3)$

$=k_1(\frac{a_{22}}{\Delta},\frac{-a_{12}}{\Delta},0)+k_2(\frac{-a_{21}}{\Delta},\frac{a_{11}}{\Delta},0)+
\lambda_3(\frac{a_{22}}{\Delta}(-a_{31})+\frac{-a_{21}}{\Delta}(-a_{32}),\frac{-a_{12}}{\Delta}(-a_{31})+\frac{a_{11}}{\Delta}(-a_{32}),1)$

$=k_1(\frac{a_{22}}{\Delta},\frac{-a_{12}}{\Delta},0)+k_2(\frac{-a_{21}}{\Delta},\frac{a_{11}}{\Delta},0)+
\lambda_3(\frac{\Delta_2}{\Delta}, -\frac{\Delta_1}{\Delta},1)$







By Corollary 3.1 the condition for $(\lambda_1,\lambda_2,\lambda_3)\in L_I(A)$ is $\lambda_3(a_{31},a_{32})\left(
                                                                                            \begin{array}{cc}
                                                                                              \frac{a_{22}}{\Delta} & \frac{-a_{12}}{\Delta} \\
                                                                                              \frac{-a_{21}}{\Delta} & \frac{a_{11}}{\Delta} \\
                                                                                            \end{array}
                                                                                          \right)=$ $\lambda_3(\frac{-\Delta_2}{\Delta},\frac{\Delta_1}{\Delta})\in \Z^2$, i.e.
$\lambda_3 \Delta_1\in \Z{\Delta}$ and $\lambda_3 \Delta_2\in \Z{\Delta}$.
Hence $\lambda_3$ is divided by $\D{\frac{\Delta}{lcm(\Delta,\Delta_1)}}$ and $\D{\frac{\Delta}{lcm(\Delta,\Delta_2)}}$. As a result $$\D{\Gamma_I\cong Z_{\frac{\Delta}{gcd(lcm(\Delta,\Delta_1),lcm(\Delta,\Delta_2))}}}.$$
\end{example}


\begin{example}
For the $(n+1)\times (n+1)$ affine Cartan matrix $$\widetilde A_n=\left(
                                                                    \begin{array}{cccccc}
                                                                      2 & -1 & 0 & \cdots &0 & -1 \\
                                                                      -1 & 2 & -1 & \cdots & 0 & 0 \\
                                                                      0 & -1 & 2 & \ddots & 0 & 0 \\
                                                                      \vdots & \vdots & \ddots & \ddots & \vdots & \vdots \\
                                                                      0 & 0 & 0 & \cdots& 2 & -1 \\
                                                                      -1 & 0 & 0 & \cdots & -1 & 2 \\
                                                                    \end{array}
                                                                  \right)$$

Let $\widetilde S=\{0,1,2,\cdots,,n-1,n\}$. Since the Dynkin diagrams are invariant under the cyclic permutation $(012\cdots (n-1) n)$, we can assume $0\not\in I$ and $1\in I$ for a subset $\emptyset \neq I\subsetneq S$. Let $k$ be the number of connected components of the Dynkin diagram of $K_I(A)$. Then $I=[i_1,j_1]\sqcup [i_2,j_2]\sqcup \cdots \sqcup [i_k,j_k]$, where $1=i_1\leq j_1< i_2\leq j_2<\cdots < i_{k}\leq j_k\leq n$.

i) If $k=1$ and $|I|=n$, then $\Gamma_I=\{e\}$.

ii) If $k=1$ and $|I|<n$, then $\Gamma_I=Z_{|I|+1}$.

iii) The case ii) can be generalized to the other case for $k>1$ and we have $$\Gamma_I\cong \Z_{j_1-i_1+2}\times \Z_{j_2-i_2+2}\times \cdots \times \Z_{j_k-i_k+2}.$$
\end{example}

The following corollary is a direct result of Theorem 2.
\begin{example}
Let $A$ be a $(n+1)\times (n+1)$ Cartan matrix of affine type and $A^\vee$ be its transposition. For each $i\in \widetilde S$, the isomorphism type of the maximal parabolic subgroup $K_{\widetilde S-\{i\}}(A)$ is $K(A_{S-\{i\}})\times S^1/\Z_{m_i}$, where $m_i$ is the integer below the $i$-th vertex of the Dynkin diagram of $A^{\vee}$ in the table of Kac\cite{Kac_82}, page 54-55.

For $i\in \widetilde S$, by deleting the $i$-th column of $A$ gives the matrix $\widetilde A_{\widetilde S-\{i\}}$. Let $\alpha_0,\alpha_1,\cdots,\alpha_n$ be the $n+1$ columns of $A^{\vee}$, then $\sum\limits_{j=0}^n m_j\alpha_j=0$. It is easy to see the $n+1$ row vectors of $A_{\widetilde S-\{i\}}$ are the transpositions of $\alpha_0,\alpha_1,\cdots,\alpha_n$ with their $i$-th components deleted. Hence we have $\sum\limits_{j=0}^n m_j e_j=0$. And $m_i e_i=-\sum\limits_{j\not=i} m_j e_j$. By Theorem 2, we obtain that $\Gamma_i\cong E_I(A)/E(A_I)\cong \Z_{m_i}$.
\end{example}

For affine Kac-Moody group of types $\widetilde A, \widetilde B, \widetilde C$ and $\widetilde D$(non-twisted or twisted cases), we always have $\det A_I\not=0$ for a subset $\emptyset \not=I\subsetneq S$. In these cases the third author designs an algorithm to compute explicitly the group $\Gamma_I$.


\end{document}